\documentclass[]{article}
\usepackage{latexsym}

\chardef\bslash=`\\ 

\newtheorem{theorem}{Theorem}[section]
\newtheorem{mytheorem}{Theorem}

\newtheorem{lemma}[theorem]{Lemma}
\newtheorem{mylemma}[mytheorem]{Lemma}

\def\Bbb{}
\def\hyp{{\Bbb H}^2}

\def\CC{\Bbb C}

\def\RR{\Bbb R}

\def\integers{\Bbb Z}

\def\tree{\bf T}

\def\Q{{\cal Q}}

\def\modud{{\cal M}_{1} (l_\delta) }
\def\modcomp{\overline{{\cal M}}_{1,1}  }
\def\modcompd{\overline{{\cal M}}_{1} (l_\delta) }
\def\teich{{\cal T}_{1,1}}
\def\teichd{{\cal T}_{1}(l_\delta)}

\def\mcg{{\cal MCG}}
\def\simp{{\cal G}_0}

\def\char{\mathrm{Hom}(\pi_1(M),\sl(\RR)) // \sl(\RR) }

\def\sech{\mathrm{sech} }
\def\arcsinh{\mathrm{arcsinh}}

\def\axis{{\mathrm axis} (A)}

\def\sys{\mathrm{sys}(x)}

\def\ZZ{\bf Z}

\def\RR{{\bf R}}

\def\sl{\mbox{SL}_2 }
\def\gl{\mbox{GL}_2 }

\def\tr{\mbox{tr\,} }

\input{epsf.sty}

\begin{document}

\title{Length series on Teichmuller space}

\author{Greg McShane}
\date{}

\maketitle

\begin{abstract}
We prove that a certain series defines a constant function using
Wolpert's  formula for the variation of the length of a geodesic
along a Fenchel Nielsen twist. Subsequently we determine the value
 viewing it as function on the the Deligne Mumford compactification
$\modcomp$ and evaluating it at the  stable curve at infinity.
\end{abstract}

\vspace{1in}

\noindent {\bf Conventions:}
\begin{enumerate}
  \item For  $\gamma$  an essential  closed
curve on a surface $l_\gamma(x)$ is the length of the geodesic
homotopic to $\gamma$ where $x$ is the point in the moduli space
determined by the  metric on the surface.
  \item For a homeomorphism  $h:M \rightarrow M$ and a geodesic
  $\gamma$,
  $h(\gamma)$ is the {\em geodesic } homotopic to the image of
  $\gamma$ under $h$.
  \item if $\gamma$ is an oriented curve then $-\gamma$ is the
  curve with the opposite orientation.
\end{enumerate}

\section{Introduction}

Let $M$ be a {\em one holed torus}. The fundamental group of $M$
is freely generated by two loops $\gamma_1,\gamma_2$ which meet in
a single point and such that their commutator is a loop, $\delta$,
around the hole. Such a surface is uniformized by a representation
$$\rho:\pi_1= \langle \gamma_1,\gamma_2\rangle \rightarrow
\sl(\RR)$$ such that the commutator of the generators is a
hyperbolic element of translation length $l_\delta$
\cite{goldmanMarkoff}. Denote by $\teichd$ the {\em Teichmuller
space of $M$} and $\modud$ the corresponding {\em moduli space}.
The mapping class group, $\mcg$, is defined to be the group of
orientation preserving diffeomorphisms up to isotopy i.e.
$\pi_0(\mbox{Diffeo}^+(M))$; it is critical  that we take only the
orientation preserving diffeomorphisms. By the work of Nielsen and
Mangel:
$$\pi_0(\mbox{Diffeo}(M)) \cong \ \mbox{Aut}(\pi_1)/\mbox{Inn}(\pi_1) \cong
\gl (\ZZ).$$ The mapping class group  is index 2 in
$\pi_0(\mbox{homeo}(M))$ and so isomorphic to $\sl(\ZZ)$.

\begin{center}
{\bf Two questions}
\end{center}

In \cite{mcshaneThesis} we showed that
$$\sum \frac{2}{1+ \exp l_\gamma} = 1$$
where the sum extends over all closed simple curves $\gamma$ on a
hyperbolic punctured torus.

\noindent {\bf Question 1:} Jorgensen  asked if the above identity
could be proved using the {\em Markoff cubic}
$$ a^2 +  b^2 + c^2
-  a b c= 0, \,\,a,b,c>2.$$ In \cite{bowditchMarkoff} Bowditch
answers this  giving a proof using a summation argument over the
edges of the tree, $\tree$, of solutions  to this equation.

\noindent {\bf Question 2:} Jorgensen also asked can the identity
be proved using Wolpert's formula for variation of length? It is
this question that we address here. By clever ''accounting"
Bowditch avoids considering the following divergent series
$$\sum_{\{a,b,c\}\in  \tree } \left( \frac{a}{bc} +\frac{b}{ca}+\frac{c}{ab}
\right).$$ The series is divergent since $T$ is infinite and,
since $a,b,c$ is a solution of the cubic, the value of each term
is $1$. Our approach is based on showing that the derivative of a
divergent series like the one above is $0$.


\begin{center}
{\bf Another  divergent series}
\end{center}

Let $\alpha, \beta$ a pair of oriented  closed simple geodesics on
$M$ meeting in exactly one point and such that the {\em signed
angle} $\alpha \vee \beta$  (see section 3) between them is
positive. By viewing $[\alpha], [\beta]$ as a basis of $H_1(M,\ZZ)
\cong \ZZ\oplus\ZZ$ \cite{mcrivin} one sees that the stabiliser of
$(\alpha, \beta)$ in $\mcg$ is trivial. Moreover, since each $g\in
\mcg$ is a homeomorphism, the pair of geodesics
$g(\alpha),g(\beta)$ again meet in a single point so $g(\alpha)
\vee g(\beta)$ is well defined and strictly positive (since $g$
preserves the orientation.)

We begin with the formal series:

$$\Q =  \sum_{g \in \mcg} g(\alpha) \vee g(\beta).$$

We determine the ``sum"  of this series using a  coset
decomposition and some hyperbolic geometry as follows.

\noindent
 {\em Step one:} we rewrite
$\Q$  as a sum over the set of coset representatives $\mcg/
\langle T_\gamma \rangle$ where $T_\gamma$ is the Dehn twist along
$\gamma$ (see section 6):

\begin{equation}\label{sum1}
  {\cal Q} = \sum_{h\in \mcg / \langle T_\gamma \rangle}
 \sum_{p\in  \langle T_\gamma \rangle} hp(\alpha) \vee hp(\beta)
= \sum_{\gamma \in \simp} \sum_{n\in \ZZ} T^n_\gamma(\gamma') \vee
T^{n+1}_\gamma(\gamma'),
\end{equation}
where the outer sum is over all oriented simple closed geodesics
$\simp$ and $\gamma'$ is any simple closed geodesic that meets
$\gamma$ exactly once. Note that, for a one holed torus, $\mcg$
acts transitively on $\simp$ and the stabiliser of $\gamma \in
\simp$ is precisely $\langle T_\gamma \rangle$ so $\mcg/ \langle
T_\gamma \rangle$ is in 1-1 correspondence with $\simp$.

\noindent {\em Step two:} we evaluate the  inner sum over $\ZZ$
using (see section 7):

\begin{mylemma}\label{lifting}
Let $\gamma$ be a simple closed curve and $\gamma'$  any simple
closed geodesic meeting $\gamma$ exactly once. Let $A\in
\rho(\pi_1)$ be an element such that $\axis/ \langle A\rangle =
\gamma$ and let $\sqrt{A}\in \sl(\RR)$ denote the square root of
$A$ . Then there exists $c\in \hyp$, and
$\hat{\alpha_n}\subset\hyp$ such that $\forall n\in \ZZ$:
\begin{enumerate}
\item
$\hat{\alpha_n}$ is a lift of $T_\gamma^n(\gamma').$
\item
if $a_n:= \hat{\alpha_n} \cap \axis $ then $a_n = (\sqrt{A})^n
(a_0).$
\item
 $c\in \hat{\alpha_n}.$

\end{enumerate}

Morever, let $\gamma^+$ (resp. $\gamma^-$ )be the geodesic passing
through $c$ and asymptotic to $\axis$ at the  attracting (resp.
repelling) fixed point of $A$. Then:
 $$ \hat{\alpha_n} \rightarrow \gamma^\pm,$$
 as $n \rightarrow \pm\infty$ and where the convergence is uniform on compact sets.

\end{mylemma}

\begin{figure}
  \centering
  \epsfbox{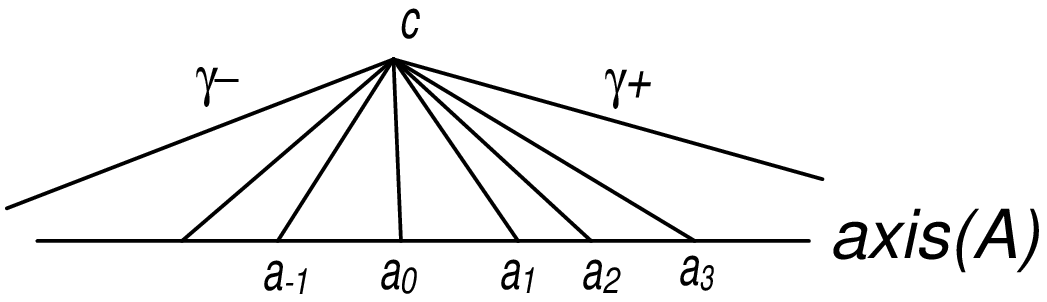}
\end{figure}
The content of the lemma is that the above diagram is a true
representation of lifts of the orbit of $\gamma'$ under
$T_\gamma$. Thus the angles in the inner sum are just the angles
between consecutive $\hat{\alpha_n}$ at the point $c$. The sum
``telescopes" over $n$ and one obtains:

\begin{equation}
\label{sum2}
  \sum_{n\in \ZZ} T^n_\gamma(\gamma') \vee
T^{n+1}_\gamma(\gamma') = \gamma^- \vee \gamma^+ =  \pi -
2\arctan\left( \frac{\cosh(l_\delta/4)}{\sinh (l_\gamma/2)}
\right) .
\end{equation}

We now determine the sum of $\Q$ using a different coset
decomposition. There is an element of order 2, $q\in \mcg$, such
that $(q(\alpha),q(\beta) ) =(\beta, -\alpha)$. Rewriting ${\cal
Q}$ as a sum over cosets of $\mcg/\langle q\rangle$, one obtains:

\begin{equation}\label{sum3}
  {\cal Q}
= \sum_{g\in \mcg/\langle q\rangle} g(\alpha) \vee g(\beta) +
gq(\alpha) \vee gq( \beta) = \sum_{\mcg/\langle q\rangle} \pi
\end{equation}
since
$$g(\alpha) \vee g(\beta) + g(\beta) \vee g(- \alpha) = \pi.$$

\noindent {\em Observation:} although this last identity clearly
implies that our series is divergent, it also suggests that the
variation of the  $\Q$ vanishes when viewed as a 1-form on
Teichmuller space. Formally one sees that:
$$ {\cal Q'} :x \mapsto  \sum_\gamma
2\arctan\left(  \frac{\cosh(l_\delta(x)/4)}{\sinh (l_\gamma(x)/2)}
\right),\,\, \teichd \rightarrow \RR $$ is constant since $d{\cal
Q'} = d{\cal Q}$.

\begin{center}
{\bf Statement of results}
\end{center}

The above illustrates a {\em formal} method for finding constant
functions defined by automorphic series over $\mcg$.  To
illustrate how this method is made rigorous we show:

\begin{mytheorem}\label{sum}
For a one holed torus $M$ :
$$\sum_\gamma \arctan\left( \frac{\cosh(l_\delta/4)}{\sinh (l_\gamma/2)}
\right) = \frac{3\pi}{2},$$ where the sum extends over all simple
closed geodesics $\gamma$ on $M$ and $l_\delta$ is the length of
the boundary geodesic $\delta$.
\end{mytheorem}

First we show that $\Q'$ is constant. From their expansions as
infinite series:
$$d{\cal Q'} = d{\cal Q} = \sum_{g} d(g(\alpha) \vee g(\beta)).$$
To conclude that the variation  vanishes one must show that the
series on the right converges  absolutely justifying the
rearrangements used above. Our point of view  is similar to that
of Kerckhoff \cite{kerckhoff} in that we do not explicitly work
with a metric, although the Weil-Petersson metric is implicit, but
with the ``Fenchel-Nielsen geometry" of the cotangent bundle. We
evaluate the pairing of $ d{\cal Q}$ with the Fenchel-Nielsen
vector field $t(\mu)$ associated to a simple closed geodesic
$\mu$. By a result of Wolpert \cite{wolpertGeo}\cite{wolpertFN}
there are finitely many simple closed geodesics $\mu_i$ such that
the associated Fenchel-Nielsen vector fields $t_{\mu_i}$ generate
the tangent space at every point in the Teichmuller space of a
surface of finite type. A 1-form vanishes iff its pairing with
these fields vanishes. Using Wolpert's formula for variation of
lengths \cite{wolpertFormula} (section 5) and elementary estimates
for the lengths of simple geodesics (section 4) we obtain as our
main theorem:

\begin{mytheorem} \label{conv}
Let $\mu$ be a simple closed geodesic $t(\mu)$ the associated
Fenchel-Nielsen vector field then the series
$$ \sum_{g\in \mcg}  d(g(\alpha) \vee g(\beta) ).t(\mu)$$
converges absolutely and its sum vanishes.
\end{mytheorem}
Absolute convergence in the usual sense for numerical series
allows one to pair terms as in (\ref{sum3}) above:
$$ d(g(\alpha) \vee g(\beta) + g(\beta) \vee g(-\alpha)  ).t(\mu) = 0,$$
and so the sum for $d\Q.t(\mu)$ vanishes identically and $\Q'$ is
constant.

Subsequently (section 8) we determine the value of the series  by
viewing it as function on the  the Deligne-Mumford
compactification $\modcompd$ and evaluating  at the  stable curve
added to obtain the compactification from $\modud$.On a
neighborhood of infinity the systole $\sys$, that is the shortest
closed geodesic, is short. To evaluate the sum we prove:

\begin{mytheorem}
\label{values} Let  $f:[0,1]\rightarrow \RR$ continuous at $1$ and
satisfying
$$f(x) = f'(0) x + O(x^{1+\delta}),$$
for some $\delta>0$.

As $\sys \rightarrow 0$,
$$
\begin{array}{lll}
\lim \sum_{\gamma} f(\sech ( l_{\gamma }/2 ) & = & f(1) +
f'(0) \lim (\sum_{\gamma'} \sech ( l_{\gamma'} /2 )) \\
&&\\
& = & f(1) + (\pi \sech(l_\delta/4))f'(0)
\end{array}
$$
where $\gamma$ varies over over all simple geodesics and $\gamma'$
over all  simple closed geodesic that meets $\sys$ exactly once.
\end{mytheorem}

\begin{center}
{\bf Generalizations}
\end{center}
There are two  generalisations of Wolpert's formula. The first is
due to Goldman \cite{goldmanSymp}\cite{goldmanHamilton} for the
representation space of a surface group into semi simple Lie group
and the second to Series \cite{seriesWolpert} for quasi-Fuchsian
deformation space. By replacing signed angle by {\em signed
complex length} \cite{seriesWolpert} one obtains a constant
function  for quasi-Fuchsian space (compare \cite{japonese}.)

Identities for higher genus surfaces  with punctures/boundary
components \cite{identities} can also be treated using this
method. In addition interesting relations can be obtained by
considering cyclic groups other than those generated by Dehn
twists.

{\bf Note.} Our approach is strikingly similar to that of Golse
and Lochak \cite{golselochak} who give an infinitesimal version of
the Selberg trace formula based on  Wolpert's formula.

\section{Markoff triples, $\teichd$}

A Markoff triple is a solution in positive integers of the {\em
Markoff cubic}:
\begin{equation}\label{markoff}
  x_1^2 +  x_2^2 + x_3^2 -  x_1  x_2 x_3 = 0,\,\, x_i >2.
\end{equation}

By work of Cohn \cite{cohn},  Haas \cite{haas} and others
\cite{seriesMarkoff} there is a correspondence between Markoff
triples and configurations $\gamma_1,\gamma_2,\gamma_3$ of simple
closed geodesics on the punctured torus  such that $\gamma_i \cap
\gamma_j, i\neq j$ is a single point. Define a {\em Markoff triple
of geodesics} to be such a configuration. Recall that a closed
geodesic, since it has no base point, only determines a conjugacy
class of the fundamental group. One can, however, view
$\gamma_1,\gamma_2$ as elements of $\pi_1(M,\gamma_1\cap\gamma_2
)$. There is a choice of orientations for $\gamma_1,\gamma_2$ them
such that:
\begin{enumerate}
  \item the commutator $[\gamma_1,\gamma_2] $ represents a loop around the puncture on
  $M$.
  \item $[\gamma_3] = [\gamma_1^{-1}\gamma_2]$ as conjugacy classes in $\pi_1$.
\end{enumerate}
\noindent{\em Observation:} The second condition is equivalent to:
$$2'.\,\,\,\,\gamma_2 = T_{\gamma_3} (\gamma_1),$$
where $T_{\gamma_3}$ is the Dehn twist round $\gamma_3$. This
point of view is important in the proof of (\ref{sum1}), see
section 6.

Now consider a representation $\rho : \pi_1(M)\rightarrow\sl(\RR)$
uniformizing a hyperbolic structure. The condition on the
commutator means that $\rho([\gamma_1,\gamma_2])$  is a parabolic
and an elementary argument shows that its trace is negative  hence
$-2$. Using the trace relations in $\sl(R)$ one computes the trace
of the commutator in terms of the traces of  $\rho([\gamma_i])$:
$$ x_1^2 +  x_2^2 + x_3^2
-  x_1  x_2 x_3 = \tr \rho([\gamma_1,\gamma_2]) + 2.$$ Thus the
triple of traces $x_i = \tr  \rho([\gamma_i])$ is always a
solution of the Markoff cubic. To obtain integer solutions one
specializes to a representation whose image is contained in
$\sl(Z)$.

A {\em generalised Markoff triple} is a solution $x_1,x_2,x_3$ of:

\begin{equation}\label{markoff2}
  x_1^2 +  x_2^2 + x_3^2 -  x_1  x_2 x_3 = -2\cosh (l_\delta/2)+
2,\,\, x_i >2.
\end{equation}

For $l_\delta = 0$, $\teichd$ is more usually denoted $\teich$; it
is a result of Keen \cite{Keen} that the solution set of
(\ref{markoff}) can be identified with the Teichmuller space
$\teich$ of $M$. Similarly for $\teichd,l_\delta > 0$ is
identified with the set of generalised Markoff triples. More
precisely, to each $x\in \teichd$ one associates a point
$[\rho_x]\in \char$, the $\sl(\RR)$ character variety of the free
group on 2 generators. The map:

$$[\rho_x] \rightarrow (\tr \rho_x ([\gamma_1]),\tr
\rho_x([\gamma_2]),\tr \rho_x ([\gamma_3]))$$ gives an embedding
of Teichmuller space into $\RR^3$ and the image satisfies
(\ref{markoff2}) ; see Goldman \cite{goldmanMarkoff} for more
details.

\section{Signed angles}

We now define the  signed angle between two geodesics at an
intersection (compare Kerckhoff \cite{kerckhoff}, Series
\cite{seriesWolpert} or Jorgensen \cite{jorg}  for a discussion of
signed complex lengths in general.) Subsequently we find an
explicit expression in terms of lengths of geodesics and study its
behaviour on Teichmuller space.

\noindent {\bf Definition:} Let $\alpha \neq \beta$ be a pair of
oriented geodesics in $\hyp$ meeting in at a point $x$. There is a
well defined {\em signed angle} between $\alpha,\beta$ at $x$;
this is a function from ordered pairs of oriented geodesics to
$]-\pi,\pi[$:
$$(\alpha,\beta) \mapsto \alpha \vee_x \beta.$$

One way to define $\alpha \vee \beta$  is to work in the disc
model for $\hyp$. After conjugating we may assume:
$$ \alpha = (-1,1),\, \beta = (-e^{i\theta},e^{i\theta}),$$
for some $\theta \in ]-\pi,\pi[$ so that $ \alpha \cap  \beta =
\{0 \}$. Now $z\mapsto e^{i\theta}z$ is the unique rotation fixing
$0$ and taking $1$ to $e^{i\theta}$ and hence $\alpha$ (oriented
in the direction from $-1$ to $1$) to $\beta$ (oriented in the
direction from $-e^{i\theta}$ to $e^{i\theta}$.) Set $\alpha
\vee_x \beta = \theta$.

For a pair of geodesics  $\alpha \neq \beta$ meeting at a point
$x$ in a surface $M$ one defines  the signed angle at $x$  by
lifting to $\hyp$. When  $\alpha \neq \beta$ meet at a single
point $x$ in the surface we shall omit $x$ and use simply $\alpha
\vee \beta$ to denote this angle.

\noindent {\bf Computation:} Let $\alpha,\beta$ be a pair of
simple closed geodesics meeting in a single point on $M$. We now
calculate the angle $\alpha \vee \beta$ in terms of $l_\alpha,
l_\beta$ and the length of the boundary $l_\delta$. Let $\gamma$
be the unique simple geodesic which satisfies $ [\gamma] +
[\alpha] =  [\beta]$ in the homology. Using the intersection form
one verifies that $\alpha,\beta,\gamma$ meet pairwise in a point
and so  form a Markoff triple. Now quotient $M$ by the elliptic
involution $J$ to obtain an orbifold $M/J$ with three cone points
one for each of the fixed points of $J$ . Simple geodesics are
invariant for this involution and the three intersection points
$\alpha\cap\beta, \beta\cap \gamma, \gamma\cap \alpha$ coincide
with the  fixed points of $J$. So the quotient of
$\alpha\cup\beta\cup\gamma$ is an embedded geodesic triangle on
$M/J$ with vertices at the 3 cone points. The sides of this
triangle have lengths $l_\alpha/2,l_\beta/2,l_\gamma/2$, and the
cosine rule for this triangle is:
$$\cosh (l_\gamma/2) = \cosh(l_\alpha/2)\cosh(l_\beta/2) -\sinh(l_\alpha/2)\sinh(l_\beta/2)\cos (\alpha \vee \beta),$$
so that
\begin{equation}\label{cosines}
 \alpha \vee \beta =\arccos\left( { \cosh(l_\alpha/2)\cosh(l_\beta/2) -\cosh (l_\gamma/2) \over
\sinh(l_\alpha/2)\sinh(l_\beta/2)} \right).
\end{equation}
One sees immediately that $\alpha \vee \beta$ is continuous on
Teichmuller space.

Finally we derive another  expression for $\alpha \vee \beta$
which will be useful in the proof of theorem \ref{conv}. Replacing
in the trace relation/Markoff cubic one obtains:
\begin{equation}\label{sines}
 \sinh^2(l_\alpha/2)\sinh^2(l_\beta/2)\sin^2 (\alpha \vee \beta )=
\cosh^2 (l_\delta/4),
\end{equation}
where $\delta$  is  the boundary geodesic.
If, $\alpha \vee \beta \in ]0,\pi/2[$ then ,

\begin{equation}\label{formula}
  \alpha \vee \beta = \arcsin \left(  {\cosh (l_\delta/4) \over
\sinh(l_\alpha/2)\sinh(l_\beta/2)}\right).
\end{equation}

\noindent {\bf Remark:} Another way of thinking of the  relation
(\ref{sines}) is as a hyperbolic version of the usual formula for
the area of a euclidean torus:

$$2 l_\alpha l_\beta \sin (\alpha \vee \beta) =
\mbox{area of torus} ,$$ where $\alpha, \beta$ are closed
Euclidean geodesics meeting in a single point at angle $\alpha
\vee \beta$.

The reader is left to check that, unfortunately, the analogous
series for the variation of the Euclidean angles does not converge
absolutely and so we obtain no new identity on the moduli space of
Euclidean structures.

\noindent {\bf Differentiability:} We now study the regularity of
$\alpha \vee \beta$ as we vary the surface over Teichmuller space.
It is well known \cite{abikoff} that for any closed geodesic
$\gamma$ the function:
$$x\mapsto l_\gamma(x), \teichd \rightarrow \RR^+$$
is differentiable and even real analytic. It is not difficult to
see from (\ref{cosines}) that $\alpha \vee \beta$ is also real
analytic.

From the expression (\ref{formula}) for the angle obtained above:

$$d (\alpha
\vee \beta)
    =  \cosh(l_\delta)/4)  \frac{\coth (l_\alpha/2) dl_\alpha  + \coth (l_\beta/2) dl_\beta  }
    { 4 (\sinh^2 (l_\alpha/2)\sinh^2 (l_\beta/2)-\cosh^2(l_\delta/4)  )
    ^{1/2}},
$$
provided $\alpha\vee\beta \neq \pi/2$ (by equation (\ref{sines}))
, that is off the subset where  $\alpha\vee\beta$ attains its
maximum. It is left to the reader that the right hand side defines
a form which extends to a continuously to the whole of Teichmuller
space by $0$ on this exceptional set.

\section{The length spectrum of simple geodesics}

We prove two lemmata used in the proof of Theorem \ref{conv} in
the next section. For a discussion of length spectra in general
see Schmutz \cite{schmutz}.

\noindent {\bf Notation:} Sections 4 and 5 deal with lengths of
geodesics and the mapping class group will not figure explicitly.
To make this clear we set
$$B^+ := \mcg.(\alpha,\beta).$$

Let $M_{g,n}$ be a hyperbolic surface of genus $g$  with $n$
punctures and $x$ the corresponding point in the moduli space. Let
$\simp$ be the set of all closed simple geodesics $\gamma \neq
\delta$. Define the {\em simple length spectrum}, denoted
$\sigma_0(x) \subset \RR^+$, to be the set $\{l_\gamma, \gamma\in
\simp\}$  counted {\em with } multiplicities. It is more useful to
think in terms of the associated {\em counting function}:
$$N(\simp,t) := \sharp \{ \gamma\in\simp, l_\gamma(x) < t\}.$$

There are two important features of the simple length spectrum:

\begin{enumerate}
  \item  The infimum over all lengths $l_\gamma(x)$ of all geodesics of
  is strictly positive and attained for a simple closed geodesic {\em the systole}, $\sys$.
  We shall also denote by  $\sys$ the length of this geodesic.

  \item $\sigma_0(M)$ is {\em discrete} that is $N(\simp,t)$ is finite for all $t\geq 0$
  and   moreover has {\em polynomial growth}
\cite{rivinCounting},\cite{Rees} that is
$$ N(\simp,t) \leq A t^{6g-6+2n},$$
for some $A = A(x)>0$.

\noindent
{\bf Remark:} If $M$ is a holed torus then it has
quadratic growth \cite{mcrivin},  $N(\simp,t) \sim At ^2$.
\end{enumerate}

\begin{lemma} \label{disc}
Let $x\in \modud$ then $\forall t >0$ there exists $N= N(t,x)>0$
such that the inequality:
$$ \sinh(l_\alpha(x)/2)\sinh(l_\beta(x)) \geq t, $$
for all but $N$ pairs $(\alpha,\beta) \in  B^+$.
\end{lemma}
{\em Proof:}
 $$ \begin{array}{lll}
    \sinh(l_\alpha) \sinh(l_\beta)
     &\geq&
           1/2(\sinh(l_\alpha/2) \sinh(\sys/2)+  \sinh(\sys/2)\sinh(l_\beta/2)) \\
     &\geq& 1/2 \sinh(\sys/2)(l_\alpha +l_\beta),

  \end{array}$$

and the lemma follows by the discreteness of the length spectrum.

\begin{lemma}
\label{growth} $\sharp \{(\alpha,\beta) \in B^+ : l_\alpha(x)
+l_\beta(x) < t \}$ grows polynomially  in $t$.
\end{lemma}

\noindent
{\em Proof:} View $B^+$ as a subset of $\simp \times
\simp$. By an elementary counting argument the growth of $N(\simp
\times \simp,t) \leq C N(\simp ,t)^2$ for some $C>0$

\noindent {\bf The Collar lemma: }Useful information about the
length spectrum can be obtained from the {\em collar lemma} see
Buser \cite{buserbook} chapter 4. Given a closed simple geodesic
$\mu$ there is an embedded collar (= regular tubular neighbourhood
of $\mu$) such that
$$\mbox{(width of collar round $\mu$)} \geq w(l_\mu), \,\,\,$$
for  $w(s) := 2\arcsinh(1/\sinh(s/2) )$. One bounds the length of
any closed geodesic $\gamma$ such that $\gamma \cap \mu \neq
\emptyset$ by theintersection number times the width the collar
round $\mu$, that is

\begin{equation}\label{collar}
 i(\gamma,\mu) \leq  \frac{l_\gamma}{w(l_\mu)},
\end{equation}
where $i(\gamma, \mu ):=\sharp \gamma \cap \mu   $ is  the {\em
geometric intersection number}.

\section{Wolperts formula and the variation of $\Q$}

Let $\mu_1,\mu_2$ be closed simple geodesics on a hyperbolic
surface. {\em Wolpert's formula} \cite{wolpertFormula} gives an
expression for the variation of $l_{\mu_1}$ along the
Fenchel-Nielsen  vector field $t(\mu_2)$  associated to $\mu_2$:
$$dl_{\mu_1}.t(\mu_2) = \sum_{z \in \mu_1 \cap \mu_2}  \cos(\theta_z),$$
where the sum is over all  the intersections between the
geodesics. An immediate corollary, also due to Wolpert
\cite{wolpert83}, is a bound  on the amplitude of the variation in
terms of intersection numbers:
$$|dl_{\mu_1}.t(\mu_2)| \leq \sum_{x \in \mu_1 \cap \mu_2} 1 =
\sharp  (\mu_1 \cap \mu_2) := i(\mu_1,\mu_2).$$  Together with the
estimates obtained in the previous section this is all that is
required to prove theorem \ref{conv}.

\noindent
{\em Proof of theorem \ref{conv}:} Fix a metric on $M$
and let $x\in \modud$ be the corresponding point in the moduli
space.

Using the formula obtained for $\alpha \vee \beta$ in section 3:
$$d (\alpha
\vee \beta)  =   \cosh(r)  \frac{ \coth (a) da  + \coth (b) db  }
    { (\sinh^2 (a)\sinh^2 (b)-\cosh(r)  ) ^{1/2}}
$$
where to simplify notation:
$$ a =l_\alpha/2, \,\, b = l_\beta/2,\,\, r =  l_\delta/4.$$

Fix a geodesic $\mu$, by Wolpert's formula:
$$
   | d (\alpha \vee \beta).t(\mu)|  \leq   \left|
   \frac{ \cosh(r) (\coth (a) i(\alpha,\mu)   + \coth (b) i(\beta,\mu) )}
   { \sinh(a)\sinh (b) (\sinh^2 (a)\sinh^2 (b)-\cosh^2(r))^{1/2}} \right|.
$$
Firstly, note that $\coth(a),\coth(b) \leq \coth(\sys/2)$ since
$a,b \geq \sys/2$. Secondly, replacing for
$i(\alpha,\mu),i(\beta,\mu)$ using (\ref{collar}) above  we obtain
the following majoration for the variation:
$$\left( \frac{\cosh(r)\coth ( \sys/2)}{ w(l_\mu)}\right).\left(\frac{l_\alpha + l_\beta}{(\sinh^2
(a)\sinh^2 (b)-\cosh^2(r)  ) ^{1/2} }\right),$$  Note that  the
leading factor does not depend on $l_\alpha, l_\beta$.

Thirdly, by the corollary to lemma  \ref{growth} for all but
finitely many pairs $(\alpha,\beta)$ in $B^+$ one has:
$$\sinh^2 (a)\sinh^2 (b)-\cosh^2(r) \geq \frac{1}{2} \sinh^2 (a)\sinh^2
(b)\geq \frac{1}{8} \exp (a+b).$$

Finally, the sum over all the configurations.
$$ \sum_{B^+}(l_\alpha + l_\beta) \exp( -1/2(l_\alpha + l_\beta) ),$$
converges since $\sharp \{(\alpha,\beta) \in B^+ : l_\alpha
+l_\beta < t \}$ grows polynomially  in $t$ by lemma \ref{growth}.

\section{Action of the mapping class group and summation of series}

We now explain Bowditch's summation argument \cite{bowditchBundle}
which decomposes a sum over $\mcg$ into a sum over orbits of all
Dehn twists.

\noindent{\bf Action on length functions:} A geodesic $\gamma$ on
$M$ determines a conjugacy class $[\gamma]$ in $\pi_1$.  A point
$x$  in Teichmuller space determines  a point in the so-called
character variety, that is an equivalence class of $\sl(\RR)$
representations of $\pi_1(M)$, $[\rho_x]$. For any representative
$\rho_x\in [\rho_x]$ the length of  $\gamma$ geodesic at $x$
satisfies:
$$2 \cosh({1\over 2}l_{\gamma}(x)) = \tr \rho_x([\gamma])$$
Now viewing  $g\in\mcg$ as a diffeomorphism of $M$, $g$ acts (on
the left) by automorphism $g_*$ on $\pi_1(M)$ and so (on the left)
on the character variety by $g:[\rho] \mapsto [\rho\circ
g^{-1}_*]$. This induces an action on the set of geodesic length
functions as follows:
$$2 \cosh (1/2 l_{g(\gamma)}(x) ) = \tr (\rho_x\circ g_*([\gamma]))
= \tr (\rho_{g^-1 x} ([\gamma]))  = 2 \cosh (1/2
l_{\gamma}(g^{-1}x) ) ,$$ for $\gamma\in\pi_1(M)$.

\noindent{\bf Summation:} Let $f:\RR^3 \rightarrow \CC$, and
$\gamma_1, \gamma_2,\gamma_3$ be a Markoff triple of geodesics.
One associates a function $\Psi: \teichd \rightarrow \CC$:
$$\Psi(x) = f(l_{\gamma_1}(x),l_{\gamma_2}(x),l_{\gamma_3}(x) ).$$
Examples of such functions are our angle $\alpha \vee \beta$,
$$(u,v,t) \mapsto \arccos\left( { \cosh(u/2)\cosh(v/2) -\cosh (t/2) \over
\sinh(u/2)\sinh(v/2)} \right),$$ and Bowditch's function,
$$B:(u,v,t)\mapsto{ \cosh(u/2)\over 2\cosh(v/2)\cosh(t/2)}.$$

We impose  a growth condition on $\Psi$ to guarantee the absolute
convergence of the sum of $\Psi$ over orbits $\mcg.x$. From the
proof of Theorem \ref{conv} a suitable condition is:
$$ | f(l_{\gamma_1}(g(x)),l_{\gamma_2}(g(x)),l_{\gamma_3}(g(x)) ) | < K \exp
-s(l_{\gamma_1}(g(x))+l_{\gamma_2}(g(x))),\,\forall g\in\mcg$$ for
some $s>0,K>0$.

Examples of such functions are $d(\alpha\vee\beta).t(\mu)$ and the
variation of Bowditch's function $d(B).t(\mu)$

By substituting in the $d(\alpha\vee\beta).t(\mu)$ in the
following theorem one obtains the  summation formula (\ref{sum1}).

\begin{lemma}
Let $\Psi,f$ be as above. Let $T $ be the Dehn twist a along
$\gamma_3$ and let $\mcg/\langle T \rangle$ denote a choice of
coset representatives for  $\langle T \rangle \subset \mcg$, then

$$
 \sum_{g\in \mcg} \Psi(g(x))  = \sum_h \sum_{n\in \ZZ}  f(l_{ hT^n( \gamma_1)}, l_{ hT^{n+1}(
\gamma_1)},l_{ h(\gamma_3)} ).$$

\end{lemma}
\noindent {\em Proof:} Consider the sum:
$$
  \begin{array}{lll}
   \sum_{g\in \mcg} \Psi(g^{-1}(x))  & = &
   \sum_{h\in \mcg/\langle T \rangle} \sum_{p\in \langle T \rangle} \Psi( (ph)^{-1}(x) )  \\
    & & \\
    & = & \sum_h \sum_{n\in \ZZ}  \Psi(h T^n (x) )\\
& & \\
& = & \sum_h \sum_{n\in \ZZ}  f(l_{ h T^n( \gamma_1)}, l_{hT^n(
\gamma_2)},l_{hT^n( \gamma_3)} )

  \end{array}
$$

Since $T$ is the Dehn twist along $\gamma_3$, $T^n(\gamma_3) =
\gamma_3$ and $T(\gamma_1) = \gamma_2$ (recall the definition of a
Markoff triple) the lemma follows. $\Box$

\section{Dehn twist orbits}

We prove the lifting lemma of the introduction. We subsequently
carry out two calculations. The first is to determine $\gamma^-
\vee \gamma^+$, thus proving equation (\ref{sum2}), in terms of
length functions and the second to determine a formula for lengths
of simple geodesics under iterated Dehn twists needed in the proof
of theorem \ref{values}.

\noindent {\em Proof of lemma \ref{lifting}:} Let $\gamma$ be a
simple closed curve and $\gamma'$  any simple closed geodesic
meeting $\gamma$ exactly once. Let $A\in \rho(\pi_1)$ be an
element such that $\axis/ \langle A\rangle = \gamma$.

For the first part of the lemma notice that since $T_\gamma^n$ is
a homeomorphism of $M$, $T_\gamma^n(\gamma')$ and $\gamma$ meet in
exactly one point. This point is necessarily one of the two
Weierstrass points on $\gamma$. Since every simple closed geodesic
passes through exactly two Weierstrass points, every curve
$T_\gamma^n(\gamma')$ passes through the unique Weierstrass point
$c$ not on $\gamma$. The point $\hat{c}\in \hyp$ is a lift of this
Weierstrass point.

Choose a lift $ \hat{\alpha_0}$ of $\gamma'$ that meets $\axis$
and choose $\hat{c}$ to be a lift of this Weierstrass point on $
\hat{\alpha_0}$ minimising the distance to $\axis$. One constructs
the geodesics $\hat{\alpha_n}$ passing through $\hat{c}$, $a_n =
(\sqrt{A})^n (a_0)$ so that these geodesics automatically satisfy
(2), (3) in the statement of the lemma. It remains to check {\em
(1)} in the statement of the lemma: that the geodesic arc joining
$\hat{c}$ to $a_n$ projects to a simple arc in the surface in the
same homotopy class as (half of) $T_\gamma^n(\gamma')$ rel the two
Weierstrass points on this latter geodesic. This is a simple
exercise left to the reader.

The second part is a simple consequence of the fact that $a_n$
converges to the attracting (resp. repelling) fixed point of $A$
as $n \rightarrow \infty$ (resp. $n \rightarrow -\infty$). $\Box$

\begin{center}
{\bf Calculation 1: angles}
\end{center}

One computes the angle $\gamma^- \vee \gamma^+$ as follows again
using the diagram \ref{myLifts}. The three geodesics $\gamma^-,
\gamma^+,\axis $ form a triangle with angles $0,0, \gamma^- \vee
\gamma^+$. Break this triangle up into two right angled triangles
formed by the perpendicular dropped from $\hat{c}$ to $\axis$ and
one half of $\axis$. One calculates $l$ the length of the
perpendicular obtaining:

\begin{equation}\label{perpendicularLength}
  \sinh(l) \sinh (l_\gamma/2) = \cosh(l_\delta/4).
\end{equation}

Now standard hyperbolic trigonometry for the right angled triangle
gives: $$\tan (\gamma^- \vee \gamma^+/2) = 1/\sinh(l).$$ One
easily obtains the quantity that appears on the right hand side of
\ref{sum2} from this.

\begin{center}
{\bf Calculation 2: lengths}
\end{center}

Note that Lemma \ref{lifting} allows one to determine the lengths
$l_{T^n(\gamma') }$ explicitly. Let $\theta$ be the signed
distance between $a_0$  and the foot of the perpendicular dropped
from $\hat{c}$ to $\axis$. Then using Pythagorus' theorem:
\begin{equation}\label{iteratedLengths}
  \cosh (l_{T^n(\gamma') }/2) = \cosh (n.l_\gamma/2 +\theta)\cosh (l).
\end{equation}

One can replace for $\cosh (l)$ from (\ref{perpendicularLength})
above.

\section{The value of a series at infinity}

We now deteremine the value  $\Q'$ by studying it in a
neighborhood of infinity in the moduli space.

Mumford and Deligne compactified the moduli space, the resulting
space is called the {\em augmented moduli space} $\modcompd$, by
adding certain {\it singular surfaces} \cite{abikoff}; these
surfaces have  double points as singularities -- each double point
is the result of  pinching an essential simple curve to a point.
Since the modular group of $M$ acts transitively on simple curves
$\neq \delta$ one adds a single point to $\modud$ to obtain the
Deligne-Mumford compactification. The Mumford-Mahler compactness
criterion says that a subset $X$ of the moduli space is precompact
iff $\sys \geq \epsilon > 0$ on $X$.  Thus a sequence of surfaces
tends to infinity in the moduli space iff the length of the
systole tends to $0$.

Our main tool is:

\begin{theorem}
\label{3 limits} Let $M$ be a punctured torus. For  $t\leq 0$, as
$\sys$ tends to $0$
$$\sum_\gamma \exp (-tl_{\gamma }) = o(\sys^{-Nt-2}),$$
where the sum extends over all simple geodesics which meet
$\alpha$ at least $N\geq 1$ times.
\end{theorem}

\noindent {\em  Proof: } We need two estimates for lengths of
curves on the torus.

Our first estimate of lengths comes via a version of the collar
lemma. Let $\alpha$ be a simple closed geodesic representing the
systole and choose $\alpha'$  so that $l_{\alpha'}$ minimizes the
lengths of geodesics $\gamma\neq\alpha$. Since
$\alpha'\neq\alpha$, $\alpha,\alpha'$ meet at least once.  The
collar lemma yields:
$$\sinh (l_{\alpha}/2)\sinh (l_{\alpha'}/2) \geq 1$$
so for $\alpha$ short one has:
$$ \exp (l_{\alpha'})  > \frac{32}{l^2_{\alpha} }.$$

For our second estimate of lengths we use the stable norm.  By
\cite{mcrivin}, for any simple closed curve $\gamma$, one has:
$$l_\gamma = \|[\gamma] \|_s, [l_\gamma]\in H_1(M,\ZZ),$$
 where $ \|.\|_s$ denotes the {\em stable
norm} on $H_1(M,\RR)$. The simple geodesics are in 1-1
correspondence with the primitive elements of $H_1(M,\ZZ)$.
Viewing $[\alpha],[\alpha']$ as a basis of the homology one writes
$ [\gamma] = m [\alpha] + n[\alpha']$ for $m,n$ coprime integers
and the intersection is equivalent to $n>N$. Thus, by dropping the
condition that the integers are coprime, one bounds the sum by:
$$2\sum_{m,  n\geq N} \exp(t \| m[\alpha] + n[\alpha']\|_s.)$$
To find upper bounds for sums of this form we need the following
inequality:

\noindent {\bf claim:} For any norm $\|.\|$  on $\RR^2$ such that
$e_1, e_2$ are the two shortest vectors with integer coefficients
there exists $K\geq 1/2$ such that:
$$\|(x,y)\| \geq K ( |x|\|e_1\| +  |y|\|e_2\|).$$

By hypothesis $[\alpha],[\alpha']$ are the shortest vectors in
$H_1(M,\ZZ)$ for the stable norm so for $K>1/2$ as above:
$$
\begin{array}{lll}
\sum_{m\geq 0,  n\geq N} \exp(t \| m[\alpha] + n [\alpha']\|_s\,)
& < &  \sum_{m,  n\geq N} \exp(t K (m\|[ \alpha ]|_s + n\|[\alpha']\|_s )\\
& & \\
& = &  (\sum_{m} \exp(t K m l_\alpha)) \,\,
(\sum_{ n\geq N} \exp(t K n l_{\alpha'})) \\
& & \\
& = & \frac{1}{1 - \exp(t K l_\alpha)}\,\,
\frac{\exp( s N K l_\alpha )}{1- \exp(t K l_{\alpha'} ) } \\
& & \\
 & < & \frac{1}{1 - \exp(t K l_\alpha)}\,\,
\frac{(R+1) l_\alpha^{-2tNK} }
 {1-l_\alpha^{-2tK} }.
\end{array}
$$
For some $R>0$. An elementary estimate shows that this latter
function  is $O(l_\alpha^{-2tNK-1})$ as $l_\alpha \rightarrow 0$.
$\Box$

\noindent {\em  Proof of theorem \ref{values}: }Let $\gamma$ be
the closed simple geodesic that realises $\sys$. The first
equation is a consequence of the preceding theorem and the fact
that the Dehn twist round $\gamma$, $T$, acts transitivley on the
set of geodesics that meet the systole exactly once. It remains to
justify:
$$\lim_{\sys\rightarrow 0} \sum_{n\in \ZZ}  \cosh(l_\delta/4)\sech ( l_{T^n(\gamma') })= \pi.$$
By a straightforward calculation using equation
(\ref{iteratedLengths}) one  has:
$$ \mbox{L.H.S.} =  \left( \sum_{n\in\integers} \sech(1/2.n.\sys + \theta)(
1/2\sys ) \right)+o(1) =  \int_{-\infty}^{\infty} \sech(u) du +
o(1) ,
$$
as $\sys\rightarrow 0$.

One evaluates the integral as usual and the theorem is proven.
$\Box$

\bibliography{myBib}
\bibliographystyle{plain}

\end{document}